\documentclass[12pt,a4paper,twoside]{article}

\newcommand{\pageformat}[6]{\setlength{\hoffset}{-1in}
                  \setlength{\voffset}{-1in}
                  \addtolength{\hoffset}{#5}
                            \addtolength{\voffset}{#6}
                            \setlength{\oddsidemargin}{#1}
                            \setlength{\evensidemargin}{#2}
                            \setlength{\textwidth}{\paperwidth}
                  \addtolength{\textwidth}{-\oddsidemargin}
                  \addtolength{\textwidth}{-\evensidemargin}
                  \addtolength{\textwidth}{-\marginparsep}
                  \addtolength{\textwidth}{-\marginparwidth}
                            \setlength{\topmargin}{#3}
                            \setlength{\textheight}{\paperheight}
                  \addtolength{\textheight}{-\topmargin}
                  \addtolength{\textheight}{-\headheight}
                  \addtolength{\textheight}{-\headsep}
                  \addtolength{\textheight}{-\footskip}
                  \addtolength{\textheight}{-#4}}
\pageformat{2cm}{3cm}{24mm}{24mm}{1pt}{0pt}

\usepackage{ifthen}
\newboolean{article}
    \setboolean{article}{true}
\newboolean{report}
\newboolean{book}
\newboolean{letter}
\newboolean{german}
\newboolean{italian}
\newboolean{nobaselinestretch}
\newboolean{nosectionappendix}
\newboolean{oldtoc}
\newboolean{nosectionequn}
\newboolean{notheorem}

\ifthenelse{\boolean{german}}{
    \usepackage{german}}{}

\usepackage[latin1]{inputenc}

\usepackage{amsmath}
\usepackage{amssymb}
\usepackage[mathscr]{eucal}

\ifthenelse{\boolean{notheorem}}{}{
    \usepackage{theorem}}

\ifthenelse{\boolean{nobaselinestretch}}{}{
    \renewcommand{\baselinestretch}{1.25}}

\newenvironment{env}[2]{\begin{#1}#2\end{#1}}{}
    \newcommand{\beq}[1]{\begin{env}{equation}{#1}}
    \newcommand{\beqn}[1]{\begin{env}{equation*}{#1}}
    \newcommand{\bal}[1]{\begin{env}{align}{#1}}
    \newcommand{\baln}[1]{\begin{env}{align*}{#1}}
    \newcommand{\bga}[1]{\begin{env}{gather}{#1}}
    \newcommand{\bgan}[1]{\begin{env}{gather*}{#1}}
    \newcommand{\bflal}[1]{\begin{env}{flalign}{#1}}
    \newcommand{\bflaln}[1]{\begin{env}{flalign*}{#1}}
    \newcommand{\bmu}[1]{\begin{env}{multline}{#1}}
    \newcommand{\bmun}[1]{\begin{env}{multline*}{#1}}
    \newcommand{\bsp}[1]{\begin{env}{split}{#1}}

    \newcommand{\eeq}{\end{env}}
    \newcommand{\eeqn}{\end{env}}
    \newcommand{\eal}{\end{env}}
    \newcommand{\ealn}{\end{env}}
    \newcommand{\ega}{\end{env}}
    \newcommand{\egan}{\end{env}}
    \newcommand{\eflal}{\end{env}}
    \newcommand{\eflaln}{\end{env}}
    \newcommand{\emu}{\end{env}}
    \newcommand{\emun}{\end{env}}
    \newcommand{\esp}{\end{env}}

\newcommand{\lf}{\vspace{2ex}}

\renewcommand{\bf}[1]{\textbf{#1}}
\renewcommand{\it}[1]{\textit{#1}}

\renewcommand{\sf}[1]{\textsf{#1}}

\renewcommand{\tt}[1]{\texttt{#1}}
\newcommand{\hl}[1]{\bf{\it{#1}}}

\newcommand{\mbf}[1]{\mathbf{#1}}
\newcommand{\msf}[1]{\text{\small$\sf{#1}$}}
\newcommand{\mtt}[1]{\mathtt{#1}}

\newcommand{\cmc}[1]{\mathcal{#1}}
\newcommand{\eus}[1]{\mathscr{#1}}
\newcommand{\euf}[1]{\mathfrak{#1}}
\newcommand{\bb}[1]{\mathbb{#1}}

\newcommand{\mfootnotesize}[1]{{\setlength{\arraycolsep}{.5ex}\text{\footnotesize$#1$}}}

\newcommand{\mtiny}[1]{{\setlength{\arraycolsep}{.3ex}\text{\tiny$#1$}}}
\newcommand{\nbd}[1]{$#1$\nobreakdash--}
\newcommand{\ol}[1]{\overline{#1}}

\newcommand{\Om}{\Omega}

\newcommand{\family}[1]{\left(#1\right)}

\newcommand{\bfam}[1]{\bigl(#1\bigr)}
\newcommand{\Bfam}[1]{\Bigl(#1\Bigr)}
\newcommand{\AB}[1]{\langle#1\rangle}
\newcommand{\bAB}[1]{\bigl\langle#1\bigr\rangle}

\newcommand{\CB}[1]{\{#1\}}
\newcommand{\bCB}[1]{\bigl\{#1\bigr\}}
\newcommand{\BCB}[1]{\Bigl\{#1\Bigr\}}
\newcommand{\SB}[1]{[#1]}
\newcommand{\bSB}[1]{\bigl[#1\bigr]}

\newcommand{\Matrix}[1]{\begin{pmatrix}#1\end{pmatrix}}

\newcommand{\fMatrix}[1]{\mfootnotesize{\Matrix{#1}}}

\newcommand{\tMatrix}[1]{\mtiny{\Matrix{#1}}}

\newcommand{\sbars}[1]{\:\bar{#1}^s\:}
\newcommand{\sodots}{\sbars{\odot}}
\newcommand{\set}[2][]{
    \ifthenelse{\equal{#1}{}}{
        \CB{#2}}{
        \CB{#1~|~#2}}}
\newcommand{\bset}[2][]{
    \ifthenelse{\equal{#1}{}}{
        \bCB{#2}}{
        \bCB{#1~|~#2}}}
\newcommand{\Bset}[2][]{
    \ifthenelse{\equal{#1}{}}{
        \BCB{#2}}{
        \BCB{#1~\big|~#2}}}
\newcommand{\zero}{\CB{0}}

\DeclareMathOperator{\ls}{\normalfont\msf{span}}
\DeclareMathOperator{\cls}{\ol{\ls}}

\DeclareMathOperator{\id}{\normalfont\msf{id}}

\newcommand{\C}{\bb{C}}

\newcommand{\N}{\bb{N}}

\newcommand{\R}{\bb{R}}

\newcommand{\cA}{\cmc{A}}
\newcommand{\cB}{\cmc{B}}
\newcommand{\cC}{\cmc{C}}
\newcommand{\cD}{\cmc{D}}

\newcommand{\cL}{\cmc{L}}

\newcommand{\sB}{\eus{B}}

\newcommand{\eH}{\euf{H}}

\newcommand{\U}{\mbf{1}}

\ifthenelse{\boolean{nosectionequn}}{}{
    \numberwithin{equation}{section}
    }

\ifthenelse{\boolean{article}\or\boolean{letter}\or\boolean{nosectionequn}}{
    \setboolean{nosectionappendix}{true}}{}
\ifthenelse{\boolean{nosectionappendix}}{}{
    \renewcommand{\appendix}{
        \chapter*{\appendixname}
        \addcontentsline{toc}{chapter}{\appendixname}
        \renewcommand{\thesection}{\Alph{section}}
        \setcounter{section}{0}}}
   
\ifthenelse{\boolean{report}\or\boolean{book}}{
    }{}

\ifthenelse{\boolean{notheorem}}{}{
        \newcommand{\mnname}{Mathematical note.}
        \newcommand{\enname}{End of the note.}
        \newcommand{\definame}{Definition.}
        \newcommand{\propname}{Proposition.}
        \newcommand{\lemname}{Lemma.}
        \newcommand{\exname}{Example.}
        \newcommand{\exername}{Exercise.}
        \newcommand{\remname}{Remark.}
        \newcommand{\obname}{Observation.}
        \newcommand{\thmname}{Theorem.}
        \newcommand{\corname}{Corollary.}
        \newcommand{\proofname}{Proof.}
    \ifthenelse{\boolean{german}}{
        \renewcommand{\mnname}{Mathematische Notiz.}
        \renewcommand{\enname}{Ende der Notiz.}
        \renewcommand{\exname}{Beispiel.}
        \renewcommand{\exername}{Übung.}
        \renewcommand{\remname}{Bemerkung.}
        \renewcommand{\obname}{Beobachtung.}
        \renewcommand{\thmname}{Satz.}
        \renewcommand{\corname}{Korollar.}
        \renewcommand{\proofname}{Beweis.}}{}
    \ifthenelse{\boolean{italian}}{
        \renewcommand{\mnname}{Nota matematica.}
        \renewcommand{\enname}{Fina della nota.}
        \renewcommand{\definame}{Definizione.}
        \renewcommand{\propname}{Proposizione.}
        \renewcommand{\exname}{Esempio.}
        \renewcommand{\exername}{Esercizio.}
        \renewcommand{\remname}{Nota.}
        \renewcommand{\obname}{Osservazione.}
        \renewcommand{\thmname}{Teorema.}
        \renewcommand{\corname}{Corollario.}
        \renewcommand{\proofname}{Dimostrazione.}

       \renewcommand{\appendixname}{Appendice}

       }{}
    \theoremheaderfont{\normalfont\bfseries}
    \theoremstyle{change}
        \theorembodyfont{\rmfamily}
            \newtheorem{emp}{}[section]
                \newcommand{\bemp}[1][]{
                    \begin{emp}\hskip-\labelsep\bf{#1}\hskip\labelsep}
                \newcommand{\eemp}{\end{emp}}
\newtheorem{itemp}[emp]{}
                \newcommand{\bitemp}[1][]{
                    \begin{itemp}\hskip-\labelsep\bf{#1}\hskip\labelsep\normalfont\itshape}
                \newcommand{\eitemp}{\end{itemp}}
            \newtheorem{mn}[emp]{\mnname}
                \newcommand{\bnm}{\begin{mn}~\begin{quotation}\renewcommand{\baselinestretch}{1}\small\noindent\ignorespaces}
                \newcommand{\enm}{\end{quotation}\hfill\bf{\enname}\end{mn}}
            \newtheorem{ex}[emp]{\exname}
                \newcommand{\bex}{\begin{ex}}
                \newcommand{\eex}{\end{ex}}
            \newtheorem{exer}[emp]{\exername}
                \newcommand{\bexer}{\begin{exer}}
                \newcommand{\eexer}{\end{exer}}
            \newtheorem{defi}[emp]{\definame}
                \newcommand{\bdefi}{\begin{defi}}
                \newcommand{\edefi}{\end{defi}}
            \newtheorem{rem}[emp]{\remname}
                \newcommand{\brem}{\begin{rem}}
                \newcommand{\erem}{\end{rem}}
            \newtheorem{ob}[emp]{\obname}
                \newcommand{\bob}{\begin{ob}}
                \newcommand{\eob}{\end{ob}}
        \theorembodyfont{\normalfont\itshape}
            \newtheorem{thm}[emp]{\thmname}
                \newcommand{\bthm}{\begin{thm}}
                \newcommand{\ethm}{\end{thm}}
            \newtheorem{prop}[emp]{\propname}
                \newcommand{\bprop}{\begin{prop}}
                \newcommand{\eprop}{\end{prop}}
            \newtheorem{cor}[emp]{\corname}
                \newcommand{\bcor}{\begin{cor}}
                \newcommand{\ecor}{\end{cor}}
            \newtheorem{lem}[emp]{\lemname}
                \newcommand{\blem}{\begin{lem}}
                \newcommand{\elem}{\end{lem}}
\newenvironment{empn}[1]{\lf\noindent\bf{#1}\ignorespaces\hskip\labelsep}{\lf}
		\newcommand{\bempn}[1]{\begin{empn}{#1}}
		\newcommand{\eempn}{\end{empn}}
		\newcommand{\bitempn}[1]{\begin{empn}{#1}\normalfont\itshape}
		\newcommand{\eitempn}{\end{empn}}
                \newcommand{\bnmn}{\begin{empn}{\mnname}~\begin{quotation}\renewcommand{\baselinestretch}{1}\small\noindent\ignorespaces}
                \newcommand{\enmn}{\end{quotation}\hfill\bf{\enname}\end{empn}}
		\newcommand{\bexn}{\begin{empn}{\exname}}
		\newcommand{\eexn}{\end{empn}}
		\newcommand{\bexern}{\begin{empn}{\exername}}
		\newcommand{\eexern}{\end{empn}}
		\newcommand{\bdefin}{\begin{empn}{\definame}}
		\newcommand{\edefin}{\end{empn}}
		\newcommand{\bremn}{\begin{empn}{\remname}}
		\newcommand{\eremn}{\end{empn}}
		\newcommand{\bobn}{\begin{empn}{\obname}}
		\newcommand{\eobn}{\end{empn}}

\newcommand{\qedsymbol}{~\rule[-0.35mm]{2mm}{2mm}}
    \newcounter{proof}[emp]
    \newenvironment{Proof}[1]{
        \vspace{1ex}
        \renewcommand{\item}[1][\stepcounter{proof}(\roman{proof})]%
            {##1\hskip\labelsep}
        \noindent\textsc{#1\hskip\labelsep}}{
        \nolinebreak\qedsymbol}
    \newcommand{\proof}[1][\proofname]{
        \begin{Proof}{#1}\ignorespaces}
    \newcommand{\qed}{\end{Proof}}
    \newcommand{\noqed}{
        \renewcommand{\qedsymbol}{}
        \end{Proof}}}
    \ifthenelse{\boolean{italian}}{
        \renewcommand{\proofname}{Dimostrazione.}}{}

\usepackage[varg]{txfonts}

\renewcommand{\thefootnote}{[\arabic{footnote}]}

\usepackage{bbm}

\begin{document}

\title{Restrictions of CP-Semigroups\\to Maximal Commutative Subalgebras\thanks{This work is supported by research funds of Italian MIUR (PRIN 2005). MS is supported by research funds of University of Molise (Dipartimento S.E.G.e S.).}\renewcommand{\thefootnote}{}\thanks{2000 AMS-Subject classification: 46L53; 46L55; 46L08; 60J25; 81S25 }}
\author{}
\author{
~\\
Franco Fagnola\\[1ex]
{\small\itshape Dipartimento di Matematica}\\
{\small\itshape Politecnico di Milano}\\
{\small\itshape Piazza Leonardo da Vinci 32}\\
{\small\itshape 20133 Milano, Italy}\\
{\small{\itshape E-mail: \tt{franco.fagnola@polimi.it}}}\\
{\small{\itshape Homepage: \tt{http://www.mate.polimi.it/$\mtt{\tilde{~}}$qp}}}\\
\\
~\\
Michael Skeide\\[1ex]
{\small\itshape Dipartimento S.E.G.e S.}\\
{\small\itshape Universit\`a\ degli Studi del Molise}\\
{\small\itshape Via de Sanctis}\\
{\small\itshape 86100 Campobasso, Italy}\\
{\small{\itshape E-mail: \tt{skeide@unimol.it}}}\\
{\small{\itshape Homepage: \tt{http://www.math.tu-cottbus.de/INSTITUT/lswas/\_skeide.html}}}\\
\\
}
\date{February 2007}

{
\renewcommand{\baselinestretch}{1}
\maketitle

\vfill



\begin{abstract}
\noindent
We give a necessary and sufficient criterion when a normal CP-map on a von Neumann algebra admits a restriction to a maximal commutative subalgebra. We apply this result to give a far reaching generalization of Rebolledo's sufficient criterion for the Lindblad generator of a Markov semigroup on $\sB(G)$.
\end{abstract}

}

\newpage


\section{Introduction}\label{intro}

The irreversible evolution of an open quantum system with associated Hilbert space $G$ is described by a (quantum) Markov semigroup on the von Neumann algebra $\cB\subset\sB(G)$ of observables. This is, from a purely mathematical point of view, a generalization of the classical Markov semigroup on some $L^\infty(\Om,\mu)$ space, that is, on a commutative von Neumann algebra $L^\infty(\Om,\mu)\subset\sB(L^2(\Om,\mu))$.

$\cB\subset\sB(G)$ may contain several commutative subalgebras. Therefore, when it turns out that one of them is invariant with respect to the action of the Markov semigroup on $\cB$, there is a classical Markov semigroup (and a classical stochastic process) embedded in the (quantum) Markov semigroup. The (classical probabilistic) information about this semigroup allows us to find valuable information on the quantum evolution. Several remarkable Markov semigroups on $\sB(G)$ admit nontrivial invariant commutative subalgebras. Indeed, all Markov semigroups arising from the stochastic limit \cite{ALV01} do. Some of them like the so-called quantum Ornstein-Uhlenbeck semigroup admit an infinite number of such invariant subalgebras; see Cipriani, Fagnola and Lindsay \cite{CFL00}.

The interest in commutative invariant subalgebras is also motivated by the study of {\em decoherence} in open quantum systems; see, for instance, Rebolledo \cite{Reb05}. This phenomenon takes place when the quantum system tends to a classical one because the off-diagonal terms (in a certain basis) of the density matrix tend to zero (and this happens on a scale faster than convergence towards an invariant state or escape to infinity). When decoherence happens the system ``chooses'' an invariant commutative algebra and the relevant evolution turns out to be given by a classical Markov semigroup.

In several important physical models on $\sB(G)$ the candidate for a commutative algebra is evident by looking at the generator. Rebolledo \cite{Reb05} gave a condition on the operators $L_i$ in the Lindblad form of the generator (see \ref{IMarkov}) for the maximal abelian algebra generated by a certain self-adjoint operator to be invariant. This is, however, only a sufficient condition. In order to determine {\em all} invariant commutative subalgebras of a given Markov semigroup, we need to find also necessary conditions. (Indeed, there are several Markov semigroups describing some phenomenological model that perhaps do not admit any nontrivial invariant commutative subalgebra. It would be good to be able to check whether, really, there is none.)

The scope of these notes is to provide a sufficient \it{and} necessary condition for that a Markov semigroup on $\cB\subset\sB(G)$ leaves invariant a commutative subalgebra. Our criterion is inspired very much by a simple generalization of Rebolledo's sufficient criterion for $\cB=\sB(G)$; see Corollary \ref{Rebcor} and Remark \ref{Rebrem}. The proof uses consequently the Hilbert module picture of the Kraus decomposition of a CP-map and the Lindblad form the generator of a Markov (or, more generally, of a CP-) semigroup. Apart from being very clear and elegant already in the case $\cB=\sB(G)$, this proof has the advantage that all statements remain true for Markov semigroups on more general von Neumann algebras $\cB\subset\sB(G)$ and commutative subalgebras $\cC\subset\cB$ that are maximal in the sense that $\cB$ does not admit bigger commutative subalgebras. However, we come always back to the case $\cB=\sB(G)$ (see Corollary \ref{Rebcor}, Remark \ref{Rebrem} and Examples \ref{Rebnonnecex} and \ref{Rebnonrepex}). We never forget that it was the sufficiency part of Corollary \ref{Rebcor} that inspired us to formulate Theorem \ref{TcharCCPthm}.

In Section \ref{prelSEC} we start with a careful introduction, explaining both the module description and how it fits together with the special versions for $\sB(G)$. In Section \ref{TSEC} we proof the criterion for a single CP-map (Theorem \ref{TcharCCPthm}). In Section \ref{LSEC} we proof the criterion for CCP-maps (Theorem \ref{LcharCCPthm}) or, what is the same, for a whole CP-semigroup.


\section{Preliminaries about von Neumann modules and GNS-con\-struc\-tions}\label{prelSEC}

\bemp[\hspace{-.6ex}.~]\label{IGNS}
Let $T\colon\cA\rightarrow\cB$ be a CP-map between unital \nbd{C^*}algebras. Since Paschke \cite{Pas73} we know how to recover $T$ in terms of a \it{GNS-construction}: Define a \nbd{\cB}valued sesquilinear map $\AB{\bullet,\bullet}$ on the vector space tensor product $\cA\otimes\cB$ by setting
\beqn{\tag{$*$}
\AB{a\otimes b,a'\otimes b'}
~:=~
b^*T(a^*a')b',
}\eeqn
turning the right \nbd{\cB}module $\cA\otimes\cB$ into a semi-Hilbert \nbd{\cB}module. The completion of the quotient  by the length-zero elements $E$ (or the strong closure in the case of von Neumann algebras) is a Hilbert (or a von Neumann) \nbd{\cB}module on which $\cA$ acts from the left by a nondegenerate representation. In other words, $E$ is a correspondence from $\cA$ to $\cB$ which we call the \hl{GNS-correspondence} associated with $T$.

The element $\U\otimes\U\in\cA\otimes\cB$ gives rise to a \hl{cyclic vector} $\xi\in E$, that is, $E=\cls\cA\xi\cB$ and we recover $T$ as $T(a)=\AB{\xi,a\xi}$. The pair $(E,\xi)$ is determined by these properties up to suitable isomorphism. We refer to $(E,\xi)$ as the \hl{GNS-construction} for $T$.
\eemp

\bemp[\hspace{-.6ex}.~]\label{IStine}
Suppose that $\cB\subset\sB(G)$ is a concrete \nbd{C^*}algebra of operators on a Hilbert space $G$. Then we may construct the Hilbert space $H:=E\odot G$ and the \hl{Stinespring representation} $\rho$ of $\cA$ on $H$ by setting $\rho(a):=a\odot\id_G$. The cyclic vector $\xi$ gives rise to a mapping $L_\xi:=\xi\odot\id_G\colon g\mapsto\xi\odot g$ in $\sB(G,H)$. We find $T(a)=L_\xi^*\rho(a)L_\xi$. This is nothing but the well-known \hl{Stinespring construction} \cite{Sti55}.
\eemp

\bemp[\hspace{-.6ex}.~]\label{IvNm}
Note that the definition of $L_\xi$ works for arbitrary elements $x\in E$. The mappings $L_x\colon g\mapsto x\odot g$ fulfill $L_x^*L_y=\AB{x,y}$ and $L_{xb}=L_xb$. We will, generally, identify $E$ as a subset of $\sB(G,H)$ by identifying $x$ and $L_x$.

If $\cB\subset\sB(G)$ is a von Neumann algebra, then also the strong closure $\ol{E}^s$ of $E$ in $\sB(G)$ is a Hilbert \nbd{\cB}module. In other words, $\ol{E}^s$ is a \hl{von Neumann \nbd{\cB}module}. If also $\cA$ is a von Neumann algebra and if $T$ is normal, then also the Stinespring representation is normal. In other words, $\ol{E}^s$ is a \hl{von Neumann correspondence} from $\cA$ to $\cB$. Von Neumann modules and correspondences (as two-sided modules) as strongly closed operator modules have been introduced in Skeide \cite{Ske00b}. The up-to-date definition is in Skeide \cite{Ske06b}. Recall that von Neumann \nbd{\cB}modules are \hl{self-dual} (that is, every bounded right linear map $E\rightarrow\cB$ has the form $x\mapsto\AB{y,x}$ for suitable $y\in E$) together with all consequences (like adjointability of all bounded module maps, existence of projections onto strongly closed submodules, and so forth).
\eemp

\bemp[\hspace{-.6ex}.~]\label{IB(G)T}
If $\cB=\sB(G)$, then $\ol{E}^s=\sB(G,H)$. ($E$ contains a norm dense subset of the finite-rank operators in $\sB(G,H)$.) If $T$ is a normal CP-map on $\sB(G)$, then the Stinespring representation is a normal nondegenerate representation of $\sB(G)$ on $\sB(H)$. The theory of these representations asserts that $H$ factors into $H=G\otimes\eH$ for some multiplicity space $\eH$ and that $\rho(a)=a\otimes\id_\eH$. In other words, $\ol{E}^s=\sB(G,G\otimes\eH)$. Let $\bfam{e_i}_{i\in I}$ denote an ONB of $\eH$. Then it is not difficult to show that the family $\bfam{\id_G\otimes e_i}_{i\in I}$ (where $\id_G\otimes e_i$ denotes the mapping $g\mapsto g\otimes e_i$) is an ONB of $\ol{E}^s$ in the obvious sense. (See \cite{Ske00b} for quasi ONBs.) Denote by $L_i:=\AB{\id_G\otimes e_i,\xi}$ the coefficients of $\xi$ with respect to this ONB. Then
\beqn{
T(b)
~=~
\sum_iL_i^*bL_i.
}\eeqn
This is a so-called \hl{Kraus decomposition} of the CP-map $T$ on $\sB(G)$.
\eemp

\bemp[\hspace{-.6ex}.~]\label{ICCP}
The formula $T_t=e^{t\cL}$ establishes a one-to-one correspondence between uniformly continuous semigroups $T=\bfam{T_t}_{t\in\R_+}$ on a the unital \nbd{C^*}algebra $\cB$ and bounded linear maps $\cL$ on $\cB$. We refer to $\cL$ as the \hl{generator} of $T$. It is well-known that $\cL$ is the generator of a \hl{CP-semigroup} (that is, all $T_t$ are CP-maps), if and only if $\cL$ is a \hl{conditionally completely positive} (CCP) map, that is, if and only if
\beqn{
\sum_{i,j}b_i^*\cL(a_i^*a_j)b_j
~\ge~
0
\text{~~~~~~whenever~~~~~~}
\sum_ia_ib_i
~=~
0,
}\eeqn
for finitely many $a_i,b_i\in\cB$. (See, for instance, Evans and Lewis \cite{EvLe77}.)
\eemp

\bemp[\hspace{-.6ex}.~]\label{IGNSCCP}
Also for CCP-maps we can construct a \hl{GNS-correspondence}. Simply take the subspace
\beqn{
(\cB\otimes\cB)_0
~:=~
\Biggl\{\sum_ia_i\otimes b_i~\Big|~n\in\N;a_i,b_i\in\cB~(i=1,\ldots,n);\sum_ia_ib_i=0\Biggr\}
}\eeqn
of $\cB\otimes\cB$, define an inner product on $(\cB\otimes\cB)_0$ by the same formula $(*)$, divide out the length-zero elements and complete as much as necessary to obtain a correspondence $E$ over $\cB$. If $\cB\subset\sB(G)$ is a concrete operator algebra, then, like for CP-maps, we may construct a Hilbert $H:=E\odot G$ with a representation $\rho(b):=b\odot\id_G$ like the Stinespring representation. This level of Stinespring-type constructions is known probably as long as generators of CP-semigroups have been studied. The GNS-type construction of the correspondence $E$, under the name \it{tangent bimodule}, is due to Sauvageot \cite{Sau89}. Sauvageot's construction is one of the very first emergencies of Hilbert modules in quantum probability. (We should like to note that Sauvageot defines the inner product on all of $\cB\otimes\cB$ by, first, projecting down to $(\cB\otimes\cB)_0$ via $a\otimes b\mapsto a\otimes b-\U\otimes ab$ and, then, defining the same inner product as above. In \cite{Sau89} this obscures slightly the origin of the left multiplication as the natural left multiplication of the \nbd{\cB}bimodule $(\cB\otimes\cB)_0$.)

As usual with generalizations of GNS-type constructions for \it{positive} structures to \it{conditionally positive} versions, finding a substitute for the \it{cyclic element} that generates everything is difficult, if not impossible. In our case, we note that the map $d\colon\cB\rightarrow E$ that sends $b\in\cB$ to the image of $b\otimes\U-\U\otimes b$ in $E$ is an \nbd{E}valued derivation on $\cB$. We also note that the range $d(\cB)$ of $d$ generates $E$ as a right Hilbert (or von Neumann) \nbd{\cB}module. In this sense, we speak of the \hl{cyclic derivation} associated with $\cL$. In the \nbd{C^*}case, in general, we do not have more than the $(E,d)$.  As usual, the pair $(E,d)$ determined by $\cL$ up to suitable isomorphism in the following sense: If $(E',d')$ is another pair such that $d'(\cB)$ generates $E'$ (in the suitable topology) as right \nbd{\cB}module and $\AB{d'(b),d'(b')}=\cL(b^*b')-\cL(b^*)b'-b^*\cL(b')+b^*\cL(\U)b'$, then $d(b)\mapsto d'(b)$ determines an isomorphism $E\rightarrow E'$ of correspondences. We refer to the pair $(E,d)$ as the \hl{GNS-construction} for $\cL$. Not that if $\cL$ is even CP, then the GNS-construction for the CCP-map $\cL$ may but need not coincide with the GNS-construction for the CP-map $\cL$.
\eemp

\bemp[\hspace{-.6ex}.~]\label{ICE}
If $\cB\subset\sB(G)$ is a von Neumann algebra, then the following result due to Christensen and Evans \cite{ChrEv79} helps a lot: A bounded derivation $d$ on $\cB$ with values in a von Neumann correspondence $E$ over $\cB$ is \hl{inner}, that is, there exists an element $\xi\in E$ such that
\beqn{
d(b)
~=~
b\xi-\xi b.
}\eeqn
(Of course, \cite{ChrEv79} do not use the language of Hilbert modules. See the appendix of Barreto, Bhat, Liebscher and Skeide \cite{BBLS04}.) $\xi$ is not unique. Replacing $E$ with $\cls^sd(\cB)\cB$ we see that $\xi$ may be chosen from the latter von Neumann submodule. Still, $\xi$ is not unique. (See Examples \ref{Rebnonnecex} and \ref{Rebnonrepex}.

Applying this crucial and hard result to the (strong closure of the) GNS-construction $(E,d)$ for a normal CCP-map $\cL$ on a von Neumann algebra $\cB\subset\sB(G)$, one obtains rather easily that $\cL$ has so-called \hl{Cristensen-Evans form}, that is,
\beqn{
\cL(b)
~=~
\cL_0(b)+b\beta+\beta^*b,
}\eeqn
where $\cL_0$ is a normal CP-map and $\beta\in\cB$. In fact, if $\xi$ is an element in $E$ such that $d(b)=b\xi-\xi b$, then $\cL_0:=\AB{\xi,\bullet\xi}$ does the job. (See the appendix of \cite{BBLS04} for a Hilbert module version of the original argument in \cite{ChrEv79}.) Like $\xi$, the Christensen-Evans form of $\cL$ is not unique.
\eemp

\bemp[\hspace{-.6ex}.~]\label{IMarkov}
$\cL$ is the generator of a \hl{Markov semigroup} $T$ (that is, $T_t(\U)=\U$ for all $t\in\R_+$), if and only if $\cL(\U)=0$. In this case, the real part $\frac{\beta+\beta^*}{2}$ of $\beta$ is necessarily given by $-\frac{\cL_0(\U)}{2}$, while the imaginary part $h=\frac{\beta-\beta^*}{2i}$ can be any self-adjoint element $h$ of $\cB$. We find that
\beqn{
\cL(b)
~=~
\AB{\xi,b\xi}-\frac{b\AB{\xi,\xi}+\AB{\xi,\xi}b}{2}+i\SB{b,h}
}\eeqn
is the sum of a purely \it{dissipative} part $\AB{\xi,b\xi}-\frac{b\AB{\xi,\xi}+\AB{\xi,\xi}b}{2}$ and a \it{hamiltonian} perturbation $i\SB{b,h}$. In the case when $\cB=\sB(G)$, so that again $E=\sB(G,G\otimes\eH)$, we find
\beqn{
\cL(b)
~=~
i\SB{b,h}+\sum_i\Bfam{L_i^*bL_i-\frac{bL_i^*L_i+L_i^*L_ib}{2}}
}\eeqn
where the $L_i=\AB{\id_G\otimes e_i,\xi}$ are the coefficients of $\xi$ with respect to some ONB $\bfam{e_i}_{i\in I}$ of $\eH$. This analogue of the Kraus decomposition of a CP-map is called the \hl{Lindblad form} of the generator $\cL$. The proof in the case $\cB=\sB(G)$ is much simpler than the general case in \cite{ChrEv79}. Lindblad's proof in \cite{Lin76} uses essentially that $\sB(G)$ may be ``approximated'' by finite-dimensional matrix algebras $M_n$, and the proof for $M_n$ uses the harmonic analysis of the (compact!) group of unitaries in $M_n$.
\eemp

\section{Globally invariant commutative subalgebras for CP-maps}\label{TSEC}

We are interested in when a normal CP-map $T$ on a von Neumann algebra $\cB\subset\sB(G)$ leaves (globally) invariant a commutative von Neumann subalgebra $\cC\ni\id_G$ of $\cB$, that is, $T(\cC)\subset\cC$.

We will give a necessary and sufficient criterion in terms of the GNS-construction $(E,\xi)$ for $T$; see \ref{IGNS}. For the proof of sufficiency we shall show that validity of our criterion implies that $\SB{T(\cC),\cC}=\zero$. In order that this suffices to show that $T(\cC)\subset\cC$, it is necessary to restrict to \hl{maximal commutative subalgebras} $\cC$ of $\cB$, in the sense that $\cC\subset\cD\subset\cB$ and $\cD$ commutative implies $\cD=\cC$. It is an easy exercise to show that this is equivalent to $\SB{b,\cC}=\zero$ $\Rightarrow$ $b\in\cC$.

We emphasize that the notion of a maximal commutative subalgebra of a von Neumann algebra $\cB$ should not be confused with the notion of a \it{maximal abelian von Neumann algebra}. A commutative von Neumann algebra on a Hilbert space $G$ is \hl{maximal abelian}, if it is a maximal commutative subalgebra of $\sB(G)$.

\bthm\label{TcharCCPthm}
Let $\cB\subset\sB(G)$ be a von Neumann algebra on the Hilbert space $G$ and let $T$ be a normal CP-map $T$ on $\cB$. Denote by $(E,\xi)$ its (strongly closed) GNS-construction. Furthermore, let $\cC\ni\id_G$ be a maximal commutative von Neumann subalgebra of $\cB$. Then $T$ leaves $\cC$ globally invariant, if and only if there exists a \nbd{*}map $\alpha\colon\cC\rightarrow\sB^a(E)$ fulfilling the following properties:
\begin{enumerate}
\item\label{Ta}
The range of $\alpha$ commutes with the left action of elements of $\cC$ on $E$, that is, for all $c_1,c_2\in\cC$ and $x\in E$ we have
\beqn{
c_1\alpha(c_2)x
~=~
\alpha(c_2)c_1x.
}\eeqn

\item\label{Tb}
For all $c\in\cC$ we have
\beqn{
\alpha(c)\xi
~=~
c\xi-\xi c.
}\eeqn
\end{enumerate}
\ethm

\proof
\it{Sufficiency}:~
If there exists a \nbd{*}map $\alpha\colon\cC\rightarrow\sB^a(E)$ fulfilling Properties \eqref{Ta} and \eqref{Tb}, then
\beq{\label{comcheck}
\SB{\AB{\xi,c_1\xi},c_2}
~=~
\AB{\alpha(c_2^*)\xi,c_1\xi}-\AB{\xi,c_1\alpha(c_2)\xi}
~=~
0
}\eeq
for all $c_1,c_2\in\cC$. As $\cC$ is a maximal commutative subalgebra of $\cB$, it follows that $\AB{\xi,c\xi}\in\cC$ for all $c\in\cC$.

\it{Necessity}:~
Suppose $T(c)\in\cC$ for all $c\in\cC$. Then the strongly closed linear subspace
\beqn{
F
~:=~
\cls^s\cC\xi\cC
}\eeqn
of $E$ is the GNS-correspondence (over $\cC$!) of $T\upharpoonright\cC$ considered as CP-map on $\cC$ with the same cyclic element $\xi$. For every $c\in\cC$ we may define the map $\delta(c)\in\sB^a(F)$
\beqn{
\delta(c)\colon
y
~\longmapsto~
cy-yc.
}\eeqn
In fact, $\delta$ is the difference of the canonical homomorphism $\cC\rightarrow\sB^a(F)$ and the map that sends $c\in\cC$ to right multiplication by $c$. The former is a \nbd{*}map into $\sB^a(F)$ and its range mutually commutes with all left actions of elements of $\cC$, because $\cC$ is commutative. The latter is a well-defined homomorphism into (actually, onto) the center of $\sB^a(F)$. So, both parts are \nbd{*}maps whose ranges commute with the left actions of elements of $\cC$. Consequently, the same is true for $\delta$.

The strongly closed linear subspace
\beqn{
F_\cB
~:=~
\cls^s\cC\xi\cB
}\eeqn
of $E$ is a von Neumann \nbd{\cB}submodule of $E$. (It is, in fact, the GNS-correspondence of $T\upharpoonright\cC$ considered as CP-map $\cC\rightarrow\cB$ with the same cyclic element $\xi$.) So, there is a projection $p\in\sB^a(E)$ onto $F_\cB$. Clearly, $F_\cB$ is invariant under the left action of $\cC$, that is, $cpx=pcpx$ for all $c\in\cC,x\in E$. From
\beqn{
cpx
~=~
pcpx
~=~
(pc^*p)^*x
~=~
(c^*p)^*x
~=~
pcx
}\eeqn
we see that $p$ commutes with the left action of all $c\in\cC$.

We note that we may identify $F_\cB$ with the tensor product $F\sodots\cB$, the von Neumann version of the tensor product over $\cC$ of the von Neumann \nbd{\cC}module $F$ with the correspondence $\cB$ from $\cC$ to $\cB$. (The left action of $\cC$ on $\cB$ is simply the restriction of the multiplication map $\cB\times\cB\rightarrow\cB$ to $\cC\times\cB$. Note that this left action is nondegenerate as $\id_G\in\cC$.) In fact, the identification $y\odot b=yb$ defines an isomorphism. Clearly, under this identification the canonical left actions of $\cC$ on $F\sodots\cB$ and on $F_\cB$ coincide.

Every element $a\in\sB^a(F)$ gives rise to an element $a\odot\id_\cB$ in $\sB^a(F\sodots\cB)=\sB^a(F_\cB)$. (On $F_\cB$ this operator acts simply as $yb\mapsto(ay)b$.) If $a$ commutes with the left action of elements of $\cC$, then so does $a\odot\id_\cB$.

Summarizing the steps we have so far in our proof of necessity, for every $c\in\cC$ we may define the operator $\alpha(c)=(\delta(c)\odot\id_\cB)p$ considered as an element in $\sB^a(E)$ that leaves $F_\cB$ invariant. As product of operators that commute with the left actions of elements of $\cC$, so does $\alpha(c)$. As $(\delta(c)\odot\id_\cB)p=p(\delta(c)\odot\id_\cB)p$ and $\delta$ is a \nbd{*}map, so is $\alpha$. Finally, since $\xi\in F\subset F_\cB$, we have
\beqn{
\alpha(c)\xi
~=~
(\delta(c)\odot\id_\cB)p\xi
~=~
(\delta(c)\odot\id_\cB)\xi
~=~
\delta(c)\xi
~=~
c\xi-\xi c.
}\eeqn
In other words, we have a \nbd{*}map $\alpha$ fulfilling Conditions \eqref{Ta} and \eqref{Tb}.\qed

\bob\label{uniob}
Note that the map $\delta$ on $F$ and its amplification to $F_\cB$ are determined uniquely by Conditions \eqref{Ta} and \eqref{Tb} restricted to $F$ and to $F_\cB$, respectively. This implies that also $\alpha$ is unique, if we put it $0$ on the complement of $F_\cB$.
\eob

\bob\label{nonGNSob}
The preceding proof does not depend on that the pair $(E,\xi)$ is the GNS-construction. It works for every vector $\xi$ in a von Neumann correspondence $E$ over $\cB$ such that $T=\AB{\xi,\bullet\xi}$.
\eob

\bcor\label{Rebcor}
Suppose that $\cB=\sB(G)$ and let $T$ be a normal CP-map on $\sB(G)$ with Kraus decomposition $T(b)=\sum_{i\in I}L_i^*bL_i$. Then $T$ leaves invariant a maximal abelian von Neumann algebra $\cC\subset\sB(G)$, if and only if for every $c\in\cC$ there exist coefficients $c_{ij}(c)\in\cC$ $(i,j\in I)$ such that
\beq{\label{cCond}
c_{ij}(c^*)
~=~
c_{ji}(c)^*
\text{~~~~~~and~~~~~~}
cL_i-L_ic
~=~
\sum_{j\in I}c_{ij}L_j.
}\eeq
\ecor

\proof
Let $E=\sB(G,G\otimes\eH)$ an arbitrary von Neumann correspondence over $\sB(G)$ and let $\xi=\sum_{i\in I}L_i\otimes e_i$ be a vector expressed with respect to some ONB $\bfam{e_i}_{i\in I}$ of $\eH$; see \ref{IB(G)T}. Consider the CP-map $T(b)=\AB{\xi,b\xi}=\sum_{i\in I}L_i^*bL_i$. Let $\alpha\colon\cC\rightarrow\sB^a(E)$ be a map and for every $c\in\cC$ define the coefficients
\beqn{
c_{ij}(c)
~:=~
\bAB{(\id_G\otimes e_i),\alpha(c)(\id_G\otimes e_j)}
}\eeqn
of $\alpha(c)$ with respect to that ONB, so that $\alpha(c)(\id_G\otimes e_j)=\sum_{i\in I}c_{ij}(c)\otimes e_i$. We observe that $\alpha(c)$ commutes with all elements of $\cC$, if and only if $c_{ij}(c)\in\cC'=\cC$. Further, $\alpha$ is a \nbd{*}map, if and only if $c_{ij}(c^*)=c_{ji}(c)^*$ for all $c\in\cC$ and $i,j\in I$. We see that there exists a \nbd{*}map $\alpha$ fulfilling Conditions \ref{TcharCCPthm}\eqref{Ta} and \eqref{Tb}, if and only if there exist $c_{ij}(c)\in\cC$ satisfying Conditions \eqref{cCond}.\qed

\brem\label{Rebrem}
The special case when $c_{ij}(c)=\delta_{ij}c_i$ for self-adjoint elements $c_i\in I$ and when $\cC$ is generated by a single self-adjoint operator $c$, is exactly Rebolledo's sufficient condition on the CP-part of a generator of a Markov semigroup in Lindblad form; see \ref{IMarkov}. In fact, it was the observation that also the Conditions \ref{cCond} are sufficient that inspired us to the present notes. But, as Corollary \ref{Rebcor} asserts, these conditions are also necessary.

We leave it as an interesting open problem, whether every suitable collection $c_{ij}(c)$ may be \it{diagonalized} by changing the ONB of $\eH$ to obtain Rebolledo's form. In the case of a general operator $\alpha(c)$ in the relative commutant of $\cC$ in $\sB(G\otimes\eH)$ this is probably not possible. But $\alpha(c)$ must also satisfy conditions with respect to the coefficients $L_i$. It is also possible that it might be necessary to consider only minimal Kraus decompositions. (See the last chapter in Parthasarathy \cite{Par92} for criteria, when the Lindblad form of a generator of a Markov semigroup on $\sB(G)$ is minimal.) In the affirmative case, this would show that Rebolledo's condition is also necessary if we allow to change the Kraus decomposition of the given CP-map.
\erem

\section{Globally invariant commutative subalgebras for CP-semi\-groups}\label{LSEC}

In this section, we are interested in when a normal uniformly continuous  CP-semigroup $T$ on a von Neumann algebra $\cB\subset\sB(G)$ leaves (globally) invariant a commutative von Neumann subalgebra $\cC\ni\id_G$ of $\cB$, that is, $T_t(\cC)\subset\cC$ for all $t\in\R_+$. It is clear that this is equivalent to $\cL(\cC)\subset\cC$ for the generator of $\cL$ of $T$.

Of course, if for some Christensen-Evans form
\beq{\label{CEf}
\cL(b)
~=~
\AB{\xi,b\xi}+b\beta+\beta^*b.
}\eeq
of $\cL$ we have that both the CP-part $\cL_0=\AB{\xi,\bullet\xi}$ and the derivation-like part $b\mapsto b\beta-\beta^*b$ leave $\cC$ globally invariant separately, then also $\cL$ leaves $\cC$ globally invariant. In particular, if
\beqn{
\cL(b)
~=~
\AB{\xi,b\xi}-\frac{b\AB{\xi,\xi}+\AB{\xi,\xi}b}{2}+i\SB{b,h}
}\eeqn
generates a Markov semigroup, then it is sufficient to check invariance for the CP-part $\cL_0$ and for the hamiltonian part $b\mapsto i\SB{b,h}$ separately; cf.\ Remark \ref{Rebrem}. But such a condition is not necessary.

\bex\label{Rebnonnecex}
Let $G=\C^2$, $\cB=\sB(G)=M_2$ and $\cC=\raisebox{.3ex}{\tMatrix{\C&0\\0&\C}}\subset\cB$. Put $L=\raisebox{.3ex}{\tMatrix{1&1\\0&1}}$ and define the CP-map $\cL_0(b)=L^*bL$ on $\cB$. By
\beqn{
\fMatrix{1&0\\1&1}\,\fMatrix{z_1&0\\0&z_2}\,\fMatrix{1&1\\0&1}
~=~
\fMatrix{z_1&z_1\\z_1&z_1+z_2}
}\eeqn
we see that $\cL_0$ does not leave $\cC$ invariant. Nevertheless, if we put $\beta=-\raisebox{.3ex}{\tMatrix{0&1\\0&0}}$, then
\beqn{
\cL_0\fMatrix{z_1&0\\0&z_2}+\fMatrix{z_1&0\\0&z_2}\beta+\beta^*\fMatrix{z_1&0\\0&z_2}
~=~
\fMatrix{z_1&0\\0&z_1+z_2},
}\eeqn
so that the CCP-map $b\mapsto L^*bL+b\beta+\beta^*b$ leaves $\cC$ globally invariant. This does not change if we normalize this map. In fact, if we put $h=\frac{\beta-\beta^*}{2i}=\frac{1}{2i}\raisebox{.3ex}{\tMatrix{0&-1\\1&0}}$, then
\beqn{
\cL(b)
~:=~
L^*bL-\frac{bL^*L+L^*Lb}{2}+i\SB{b,h}
}\eeqn
generates a Markov semigroup on $\cB$.
\bmun{
\cL\fMatrix{z_1&0\\0&z_2}
~=~
\fMatrix{z_1&z_1\\z_1&z_1+z_2}-\frac{1}{2}\family{\fMatrix{z_1&0\\0&z_2}\,\fMatrix{1&1\\1&2}+\fMatrix{1&1\\1&2}\,\fMatrix{z_1&0\\0&z_2}}+\frac{1}{2}\family{\fMatrix{z_1&0\\0&z_2}\,\fMatrix{0&-1\\1&0}-\fMatrix{0&-1\\1&0}\,\fMatrix{z_1&0\\0&z_2}}
\\
~=~
\fMatrix{z_1&z_1\\z_1&z_1+z_2}-\frac{1}{2}\family{\fMatrix{z_1&z_1\\z_2&2z_2}+\fMatrix{z_1&z_2\\z_1&2z_2}}+\frac{1}{2}\family{\fMatrix{0&-z_1\\z_2&0}-\fMatrix{0&-z_2\\z_1&0}}
~=~
\fMatrix{0&0\\0&z_1-z_2}
}\emun
shows that the restriction of $\cL$ to $\cC$ generates a classical \it{two-state death process}, although neither the CP-part $\cL_0$ nor the hamiltonian part $b\mapsto i\SB{b,h}$ leave invariant $\cC$, separately.
\eex

We want to give a sufficient \it{and} necessary condition.

\bthm\label{LcharCCPthm}
Let $\cB\subset\sB(G)$ be a von Neumann algebra on the Hilbert space $G$ and let $\cL$ be a (bounded) normal CCP-map $\cL$ on $\cB$. Denote by $(E,d)$ its (strongly closed) GNS-construction. Furthermore, let $\cC\ni\id_G$ be a maximal commutative von Neumann subalgebra of $\cB$. Then $\cL$ leaves $\cC$ globally invariant, if and only if there exist an element $\zeta\in E$ that reproduces $d\upharpoonright\cC$ as
\beqn{
d(c)
~=~
c\zeta-\zeta c,
}\eeqn
a \nbd{*}map  $\alpha\colon\cC\rightarrow\sB^a(E)$ and a self-adjoint element $\gamma\in\cC$ such that the following conditions are satisfied:
\begin{enumerate}
\item\label{L2a}
The range of $\alpha$ commutes with the left action of elements of $\cC$ on $E$, that is, for all $c_1,c_2\in\cC$ and $x\in E$ we have
\beqn{
c_1\alpha(c_2)x
~=~
\alpha(c_2)c_1x.
}\eeqn

\item\label{L2b}
For all $c\in\cC$ we have
\beqn{
\alpha(c)\zeta
~=~
c\zeta-\zeta c.
}\eeqn

\item\label{L1}
For all $c\in\cC$ we have
\beqn{
\cL(c)-\AB{\zeta,c\zeta}
~=~
\gamma c.
}\eeqn
\end{enumerate}
\ethm

\proof
\it{Sufficieny}:~
By Theorem \ref{TcharCCPthm}, Conditions \eqref{L2a} and \eqref{L2b} imply that the CP-map $b\mapsto\AB{\zeta,b\zeta}$ on $\cB$ leaves $\cC$ globally invariant. By Condition \eqref{L1}, the same is true for $\cL$.

\it{Necessity}:~
Suppose $\cL(c)\in\cC$ for all $c\in\cC$. Then $F:=\cls^sd(\cC)\cC\subset E$ is just the GNS-correspondence of $\cL\upharpoonright\cC$ considered as CCP-map on $\cC$; see \ref{IGNSCCP}. By \cite{ChrEv79} there exist $\zeta\in F$ and $\gamma_0\in\cC$ such that $d(c)=c\zeta-\zeta c$ and $\cL(c)=\AB{\zeta,c\zeta}+c\gamma_0+\gamma_0^*c$; see \ref{ICE}. By commutativity of $\cC$, we have $c\gamma_0+\gamma_0^*c=\gamma c$ with $\gamma=\gamma+\gamma^*$. This shows Condition \eqref{L1}.

As in the proof of necessity in Theorem \ref{TcharCCPthm}, by setting
\beqn{
\delta(c)\colon
y
~\longmapsto~
cy-yc.
}\eeqn
we define a \nbd{*}map $\delta$ from $\cC$ into the \nbd{\cC}bilinear operators on $F$, that fulfills $\delta(c)\zeta=c\zeta-\zeta c$. Again, by $p\in\sB^a(E)$ we denote the projection onto the von Neumann \nbd{\cB}submodule
\beqn{
F_\cB
~:=~
\cls^sd(\cC)\cB
~=~
\ol{\,\BCB{(c\zeta-\zeta c)b\colon c\in\cC,b\in\cB}\,}^{\,s}
}\eeqn
of $E$ generated by $F$. By
\beqn{
c_1(c_2\zeta-\zeta c_2)b
~=~
(c_1c_2\zeta-\zeta c_1c_2)b-(c_1\zeta-\zeta c_1)c_2b
}\eeqn
we see that $F_\cB$ is invariant under the left action of $\cC$ so that, once more, $p$ commutes with the left action of all $c\in\cC$. Also here, we may identify $F_\cB$ with the tensor product $F\sodots\cB$. In conclusion, the map $\alpha$ defined by setting $\alpha(c)=(\delta(c)\odot\id_\cB)p$ considered as an element in $\sB^a(E)$ fulfills Conditions \eqref{L2a} and \eqref{L2b}.\qed

\lf
It is noteworthy that this condition does not involve any Christensen-Evans form for $\cL$ but only a Christensen-Evans form for $\cL\upharpoonright\cC$, if the latter exists. Even if we know a Christensen-Evans form \eqref{CEf} for $\cL$, this does not really help to apply Theorem \ref{LcharCCPthm}. In Theorem \ref{TcharCCPthm}, $\xi$ is given from the beginning and it is essentially unique. Also, in Theorem \ref{TcharCCPthm} there is not much choice how to define $\alpha$. Here, before we can try to find $\alpha$, we must first find a candidate for $\zeta$. We know that, if it exists, then we can find one in $\cls^s d(\cC)\cC$.

\bex\label{Rebnonrepex}
Let us return to Example \ref{Rebnonnecex}. We easily verify that $E=\cB$ and $\xi=L$. Definitely, $\xi$ cannot serve as $\zeta$, because $\AB{\xi,\xi}=L^*L=\raisebox{.3ex}{\tMatrix{1&1\\1&2}}\notin\cC$. We observe that
\beqn{
d\fMatrix{z_1&0\\0&z_2}
~=~
\fMatrix{z_1&0\\0&z_2}\,\fMatrix{1&1\\0&1}-\fMatrix{1&1\\0&1}\,\fMatrix{z_1&0\\0&z_2}
~=~
\fMatrix{z_1&z_1\\0&z_2}-\fMatrix{z_1&z_2\\0&z_2}
~=~
\fMatrix{0&z_1-z_2\\0&0},
}\eeqn
so that $d(\cC)\cC=\C\,\raisebox{.3ex}{\tMatrix{0&1\\0&0}}$. If we put $L'=\raisebox{.3ex}{\tMatrix{0&1\\0&0}}$, we see that $L'^*\raisebox{.3ex}{\tMatrix{z_1&0\\0&z_2}}L'=\raisebox{.3ex}{\tMatrix{0&0\\0&z_1}}$ and
\beqn{
L'^*\fMatrix{z_1&0\\0&z_2}L-\fMatrix{0&0\\0&1}\,\fMatrix{z_1&0\\0&z_2}
~=~
\fMatrix{0&0\\0&z_1-z_2}
}\eeqn
gives back $\cL\upharpoonright\cC$ from Example \ref{Rebnonnecex}. This shows that we may put $\zeta=L'$ and $\gamma=\raisebox{.3ex}{\tMatrix{0&0\\0&1}}$.

Note that $L=L'+\id_G$. So, $L^*bL=(L'^*+\id_G)b(L'+\id_G)=L'^*bL'+bL'+L'^*b+b$ and, therefore,
\beqn{
\cL(b)
~=~
L'^*bL'+bL'+L'^*b+b-\frac{bL^*L+L^*Lb}{2}+i\SB{b,h}.
}\eeqn
The real part of $L'$ is $\frac{1}{2}\raisebox{.3ex}{\tMatrix{0&1\\1&0}}$, the imaginary part $\frac{1}{2i}\raisebox{.3ex}{\tMatrix{0&1\\-1&0}}=-h$. We find
\bmun{
\cL(b)
~=~
L'^*bL'+b-\frac{1}{2}\family{b\fMatrix{1&0\\0&2}+\fMatrix{1&0\\0&2}b}
~=~
L'^*bL'-\frac{1}{2}\family{b\fMatrix{0&0\\0&1}+\fMatrix{0&0\\0&1}b}
\\
~=~
L'^*bL'-\frac{bL'^*L'+L'^*L'b}{2}.
}\emun
\eex

The example tells us two things: Firstly, it may happen that $\zeta$ can replace $\xi$. That is, not only $d(c)=c\zeta-\zeta c$ for all $c\in\cC$ but even $d(b)=b\zeta-\zeta b$ for all $b\in\cB$. Secondly, an inconvenient choice for $\xi$ may even cause a hamiltonian part that, otherwise, would not be there.

For the first observation, it would certainly be good, if we could proof the converse, namely, for every $\cL$ leaving invariant a maximal commutative subalgebra $\cC$ of $\cB$ there is a Christensen-Evans form such that the CP-part alone leaves $\cC$ invariant. Presently, we do not yet have a feeling whether the answer might be affirmative.

For the second observation, this is settled by the following probably well-known lemma: Any nontrivial hamiltonian part in a Christensen-Evans form of $\cL$ must be there to compensate missing invariance of the CP-part $\cL_0$.

\blem
Suppose $\cC\subset\sB(G)$ is a commutative von Neumann algebra, and let $h\in\sB(G)$ be such that $\SB{c,h}\in\cC'$ for all $c\in\cC$. Then $\SB{c,h}=0$. In particular, if $\cC$ is a maximal commutative subalgebra of the von Neumann algebra $\cB\subset\sB(G)$ and $h\in\cB$, then $h\in\cC$.
\elem

\proof
Let $p,q\in\cC$ denote projections such that $pq=0$. Then
\beqn{
0
~=~
\bSB{p,\SB{q,h}}
~=~
pqh-phq-qhp+hqp
~=~
-phq-qhp.
}\eeqn
Multiplying with $p$ from one side, we find $phq=0=qhp$. It follows $q\SB{p,h}=0=\SB{p,h}q$ and, in particular, $(1-p)\SB{p,h}=0=\SB{p,h}(1-p)$. Further we compute
\baln{
\SB{p,h}
~=~
p\SB{p,h}
&
~=~
ph-php,
&
\SB{p,h}
~=~
\SB{p,h}p
&
~=~
php-hp.
}\eeqn
Adding the two equations we find
\beqn{
2\SB{p,h}
~=~
ph-hp
~=~
\SB{p,h}.
}\eeqn
In other words, $\SB{p,h}=0$ for every projection $p\in\cC$. As every $c\in\cC$ is the norm limit of linear combinations of projections in $\cC$, it follows $\SB{c,h}=0$ for all $c\in\cC$.\qed

\lf
It appears appealing to check our construction against the complete class of Markov semigroups on $M_2$ that leave invariant the diagonal subalgebra and the off-diagonal subspace and that admit an invariant state, as determined explicitly by Carbone \cite{Car04}.


\setlength{\baselineskip}{2.5ex}


\newcommand{\Swap}[2]{#2#1}\newcommand{\Sort}[1]{}
\providecommand{\bysame}{\leavevmode\hbox to3em{\hrulefill}\thinspace}
\providecommand{\MR}{\relax\ifhmode\unskip\space\fi MR }
\providecommand{\MRhref}[2]{%
  \href{http://www.ams.org/mathscinet-getitem?mr=#1}{#2}
}
\providecommand{\href}[2]{#2}


\end{document}